\documentclass[final]{siamltex}
\usepackage{epsfig}

\def\Id{\mbox{\boldmath$I$}}

\def\R{\mbox{\boldmath$R$}}
\def\E{\mbox{\boldmath$E$}}
\def\A{\mbox{\boldmath$A$}}
\def\K{\mbox{\boldmath$K$}}
\def\P{\mbox{\boldmath$P$}}
\def\Q{\mbox{\boldmath$Q$}}

\def\M{\mbox{\boldmath$M$}}
\def\W{\mbox{\boldmath$W$}}
\def\H{\mbox{\boldmath$H$}}

\def\D{\mbox{\boldmath$D$}}
\def\U{\mbox{\boldmath$U$}}
\def\L{\mbox{\boldmath$L$}}

\def\Minv{\mbox{\boldmath$M$}^{-1}}
\def\Ainv{\mbox{\boldmath$A$}^{-1}}
\def\Pinv{\mbox{\boldmath$P$}^{-1}}

\newcommand{\ie}{{\em i.e.,}}

\title{
	Multi-scale linear solvers for very large systems
	derived from PDEs%
	\footnote{This work presented at the
		Fifth Copper Mountain conference on Iterative Methods,
		April, 1998.
		}
}

\author{
	Klaus Lackner%
		\thanks{Theoretical Division, Mail Stop B-216,
			Los Alamos National Laboratory,
		 	Los Alamos, NM 87544 (\texttt{ksl@lanl.gov}).}
        \and Ralph Menikoff%
		\thanks{Theoretical Division, Mail Stop B-214,
			Los Alamos National Laboratory,
		 	Los Alamos, NM 87544 (\texttt{rtm@lanl.gov}).}
}

\begin{document}

\maketitle

\begin{abstract}
We present a novel linear solver that works well for large systems
obtained from discretizing PDEs.
It is robust and, for the examples we studied, the computational effort scales
linearly with the number of equations.
The algorithm is based on a wavelength decomposition
that combines conjugate gradient,
multi-scaling and iterative splitting methods into a single approach.
On the surface, the algorithm is a simple preconditioned conjugate gradient
with all the sophistication of the algorithm in the choice of the
preconditioning matrix.
The preconditioner is a very good approximate inverse of the  linear operator.
It is constructed from the inverse of the coarse grained linear operator
and from smoothing operators that are based on an operator splitting
on the fine grid.  The coarse graining captures the long wavelength
behavior of the inverse operator while the smoothing operator captures
the short wavelength behavior.  The conjugate gradient iteration
accounts for the coupling between long and short wavelengths.
The coarse grained operator corresponds to a lower resolution approximation
to the PDEs.  While the coarse grained inverse is not known explicitly,
the algorithm only requires that the preconditioner can be a applied to
a vector.  The coarse inverse applied to a vector can be obtained
as the solution of another preconditioned conjugate gradient solver
that applies the same algorithm to the smaller problem.
Thus, the method is naturally recursive.
The recursion ends when the matrix is sufficiently small
for a solution to be obtained efficiently with a standard solver.
The local feedback provided by the conjugate gradient step at every level
makes the algorithm very robust.
In spite of the effort required for the coarse inverse,
the algorithm is efficient because the increased quality
of the approximate inverse greatly reduces the number of times
the preconditioner needs to be evaluated.
A feature of the algorithm is that the transition between coarse grids
is determined dynamically by the accuracy requirement of the conjugate
gradient solver at each level.  Typically, later iterations on the
finer scales need fewer iterations on the coarser scales and the
computational effort is proportional to $N$ rather than $N \log N$,
where $N$ is the number of equations.
We have tested our solver on the porous flow equation.
On a workstation we have solved problems on
grids ranging in dimension over 3 orders of magnitude,
from $10^3$ to $10^6$, and found that the linear scaling holds. 
The algorithm works well, even when the
permeability tensor has spatial variations exceeding a factor of $10^9$.
\end{abstract}



\pagestyle{myheadings}
\thispagestyle{plain}
\markboth{K. LACKNER AND R. MENIKOFF}{Multi-Scale Linear Solvers}

\section{Introduction}
We present a novel linear solver that works well for large systems
obtained from discretizing PDEs.
It is robust and for the examples we studied, the computational effort scales
linearly with the number of equations.
The algorithm is based on a wavelength decomposition
which combines conjugate gradient,
multiscaling and iterative splitting methods into a single approach.
On the surface, the algorithm is a simple preconditioned conjugate gradient
with all the sophistication of the algorithm in the novel choice of the
preconditioning matrix.
The preconditioner is an approximate inverse of the  linear operator
constructed from the inverse of the coarse grained linear operator
and a smoothing operator based on an iterative expansion of the
operator on the fine grid.  To transform the coarse grained operator
to the fine grid it is pre- and post-multiplied by reduction and prolongation
operators. The coarse graining captures the long wavelength behavior
of the inverse operator
while the smoothing operator captures the short wavelength behavior.
The conjugate gradient iteration accounts for the coupling between long
and short wavelengths.

The coarse grained operator corresponds to a lower resolution approximation
to the PDEs.  While its inverse is not known explicitly,  the conjugate
gradient algorithm only requires that the preconditioner can be applied to
a vector.  The inverse is obtained
as the solution of another preconditioned conjugate gradient solver
that applies the same algorithm to the smaller problem.
Thus, the method is naturally recursive.
The recursion ends when the matrix is sufficiently small
for a solution to be obtained efficiently with a standard solver.

Basing a preconditioner on the inverse of a simpler problem is
akin to using an incomplete factorization.
As with an incomplete factorization, our preconditioner only requires
knowledge of sparse matrices.
But rather than a simple forward and backward substitution used to evaluate
the inverse of the product of a lower and upper triangular matrix,
our preconditioner requires an iterative method.
Despite the extra effort required for the coarse inverse,
the algorithm is efficient because the increased quality
of the approximate inverse greatly reduces the number of times
the preconditioner needs to be evaluated.
In the cases we studied the computational effort is dominated
by the operations performed on the finest scale.

As with the multigrid method, by utilizing several scales we obtain rapid
convergence of long wavelengths.  In contrast to the way multigrid algorithms
are typically applied,
the recursive use of the conjugate gradient algorithm at each stage enforces an
accurate transition from one scale to the next.
In particular, the algorithm naturally accounts for the coupling
between the long and short wavelengths introduced at every level of refinement.
Thus, it avoids the accumulation of errors that stem from the interaction
terms between the different levels.
In addition, the transition between levels is determined dynamically
rather than preprogrammed with a V-cycle or W-cycle.
The algorithm is robust and works well without
fine tuning even when the PDE is far from diagonal in Fourier space.

\section{Algorithm}
Let the linear system be given by
\begin{equation}
	\A x = b	\;.	\label{eq:2.1}
\end{equation}
With a preconditioning matrix $\Minv$, the system to be solved is
\begin{equation}
	(\Minv \A) x = \Minv b	\;.	\label{eq:2.2}
\end{equation}
Our aim is to determine a good approximate inverse of $\A$
that can be used for $\Minv$.
This reduces the condition number of the system
and hence enhances the convergence rate of any iterative solver.
We use a conjugate gradient algorithm
since it has the advantage that convergence can greatly be enhanced
even if the preconditioner fails to account for a small number of modes.

To develop an approximate inverse,
we start with a standard operator splitting of the form,
see for example \cite{Axelsson},
\begin{equation}
	\A = \P - \Q \;.					\label{eq:2.3}
\end{equation}
We assume that $A$ and $\P$ are positive and symmetric.
For our algorithm, it is necessary that the inverse of $\P$ applied to a vector
can be evaluated efficiently and that $\Pinv  \Q$ damps out short wave lengths.
The inverse is given formally by
\begin{equation}
	\Ainv  = (\Id - \Pinv \Q)^{-1} \Pinv 
	       = \sum_{n=0}^\infty (\Pinv \Q)^n \Pinv \;. \label{eq:2.4}
\end{equation}
The series converges when $|| \Pinv \Q || < 1$, where $|| \cdot ||$
denotes the $L^2$ matrix norm.  Frequently, the matrix arising from
discretizing a PDE is diagonally dominated and it can be shown that
standard splittings, such as Jacobi or Gauss-Seidel, do indeed satisfy
the sufficient condition for the series to converge.
However, as the resolution of the discretization increases $|| \Pinv \Q ||$
typically approaches 1. Then, the convergence of the series for $\Ainv$
is very slow.

We observe the following identity
\begin{equation}
	\Ainv  = (\Pinv \Q)^m \Ainv  (\Q\Pinv )^m
			+ \sum_{k=0}^{2m-1}(\Pinv \Q)^k \Pinv 
							\label{eq:2.5}
\end{equation}
for any integer $m \ge 1$.  This identity
is easily proved by substituting (\ref{eq:2.4})
for $\Ainv $ on the right hand side of the equation
and resumming the series.
If one replaces $\Ainv$ on the right hand side of the equation with
a matrix $\W$ corresponding to an approximation of $\Ainv$,
then one obtains an approximate inverse better than $\W$:
\begin{equation}
	\Minv  = (\Pinv \Q)^m \W (Q\Pinv )^m
			+ \sum_{k=0}^{2m-1}(\Pinv \Q)^k \Pinv  \;.
							\label{eq:2.6}
\end{equation}
A measure of the improvement of the inverse can be obtained from the formula
\begin{equation}
	\A \Minv  -\Id = (\Q \Pinv )^m \cdot(\A \W-\Id)\cdot(\Q\Pinv )^m
						\; .\label{eq:2.7a}
\end{equation}
Consequently, $||\A \Minv  -\Id|| \le || \Pinv \Q||^{2m}\, || \A \W - \Id||$.
Thus, with respect to this norm, $A \Minv $ is closer to the identity than
$\A\W$ provided that $||\Pinv \Q|| < 1$.  For large systems
$||\Pinv \Q||$ is near 1 and for $\Minv $ to be an effective preconditioner
$\W$ and $\Pinv \Q$ need to be complementary in the sense
that they approximate different parts of the spectrum of $\A$.
In this case,
$\Minv$ can be a much better inverse of $\A$ than
either $\W$ or the series given in (\ref{eq:2.4}) truncated after $2m$ terms.

Equation (\ref{eq:2.6}) forms the basis for a number of approximation schemes.
Tatebe \cite{Tatebe} pointed out that $\Minv $ will be symmetric
positive definite if $\W$, $\P$ , and $\Q$  are all symmetric and
positive.  We note that $\W$ need be only positive and not
positive definite for $\Minv $ to be a valid preconditioner for a
conjugate gradient algorithm.  Consequently, the coarse grained inverse
can be used for $\W$ even though its null space is non-empty.

\subsection{Polynomial Preconditioner}
The simplest choice is $\W = \Pinv$. In this case, $\Minv$ corresponds to a 
conventional polynomial preconditioner; \ie the first $2m+1$ terms of the
series for $\Ainv$ in (\ref{eq:2.4}).
At high resolution, $||\Pinv  \Q||$ is typically
close to 1 and $\Minv$ is a poor preconditioner.
The underlying reason is that $\P$  connects only neighboring grid points
and consequently $\Minv$ provides a poor approximation of $\Ainv$
at long wave lengths.  This is born out by experience showing
that the number of iterations grows rapidly with the dimension of the
system.

\subsection{Multi-Grid Preconditioner}
A better choice is that advocated by Tatebe \cite{Tatebe} which aims to
account for the long wavelengths by basing the preconditioner
on a single step of a multigrid algorithm. Let $G_k$ be a sequence of
successively coarser grids, $\R_{k+1,k}$ and $\E_{k,k+1}$ be the reduction
and prolongation operators connecting the grids $G_k$ and $G_{k+1}$,
$\A_k$ the coarsened operator on the grid $G_k$,
and $\P_k$ and $\Q_k$ the splitting of $\A_k$.
Here $k=0$ corresponds to the finest mesh and $k_{\rm max}$
to the coarsest mesh.  The preconditioner can be defined
recursively as
\begin{equation}
	\M_k^{-1} = (\H_k)^{m} (\E_{k,k+1} \M_{k+1}^{-1} \R_{k+1,k}) (\H_k^T)^m
			+ \sum_{j=0}^{2m-1}(\H_k)^j \P_k^{-1} \;,
							\label{eq:2.9}
\end{equation}
where $\H_k = \P_k^{-1}\Q_k$.  Provided that the reduction and prolongation
operators are chosen such that $\R_{k+1,k}^T = \E_{k,k+1}$, the approximate
inverse at every level, $\W_k = \E_{k,k+1} \M_{k+1}^{-1} \R_{k+1,k}$
is symmetric.
On the coarsest mesh, the problem can be solved exactly.  Thus, the recursion
ends with $\M_{k_{\rm max}}^{-1} = \A_{k_{\rm max}}^{-1}$.
The standard multigrid algorithm uses $\M_0^{-1}$ as an approximate
inverse in conjunction with a simple iterative solver.
By contrast Tatebe's algorithm uses $\M_0^{-1}$ as a preconditioner
for a conjugate gradient solver.

For $\M_0^{-1}$ to be a good approximate inverse of $\A$,
the splitting must be chosen
such that $\H_k$ smoothes the shorter wavelengths to complement $\W_k$
which operates only on the longer wavelengths.  Thus, the underlying
rationale for the multigrid preconditioner to be a good approximate
inverse is a wavelength decomposition.  When the system is large enough
to require many coarsening levels, two problems can arise.  The first
is that many iterations are needed to account for the coupling between
short and long wavelengths.  The second is that errors due to the
coupling can accumulate at each level to such an extent that the
iterative improvement to the solution saturates and the
the desired tolerance can not be achieved.

\subsection{Recursive Multi-Scale Conjugate Gradient}
The difficulties with the multi-grid preconditioner can be overcome by
defining the preconditioner in term of the exact inverse of the
coarsened operator $A_c$ as follows
\begin{equation}
	\Minv  = (\Pinv \Q)^{m} (\E \A_c^{-1}\R) (\Q \Pinv )^m
			+ \sum_{j=0}^{2m-1}(\Pinv \Q)^j \Pinv  \;,
							\label{eq:2.10}
\end{equation}
In effect, $\W = \E \A_c^{-1} \R$ and the evaluation of $\A_c^{-1}$ on
a vector is computed by applying the same algorithm to the coarse
grained operator.  Again the recursion ends by solving the problem
exactly on the coarsest level.  Using the exact inverse $\A_c^{-1}$ for
$\W$ is the best preconditioner based on a single level of scaling.

An important feature of this preconditioner is that the long and short
wavelengths are coupled at each level by a conjugate gradient
solver before proceeding to the next finer grid.  This prevents truncation
errors from the coarsening of the operator at each level from accumulating.
In particular, if the coarsened operator does poorly on a few modes,
which typically occurs when the coefficients of the underlying PDE are
discontinuous, then these modes will be corrected by
the conjugate gradient solver with only a small penalty.

The efficiency of the algorithm is due to the high quality of the
approximate inverse.  This results in a preconditioned matrix
$\Minv \A$ with a low condition number and consequently a small number
of iterations for the conjugate solvers at every level. The condition
number for positive symmetric matrices is the ratio of the largest to
smallest eigenvalues.  The coarsening increases the smaller eigenvalues
associated with the long wavelengths and the iterative splitting
decreases the larger eigenvalues associated with short wavelengths.
Squeezing the eigenvalues together greatly decreases the condition
number of the preconditioned matrix.
Of course the reduced number of iterations must overcome any
increase in the cost per iteration.  The numerical examples below show that
this is indeed the case.

The efficiency of the conjugate gradient algorithm comes from the simple
recursion relations for the conjugate directions.  The orthogonality of
the conjugate directions is a consequence of the linearity of the
preconditioner.  A potential draw back of a preconditioner that depends on
a linear solver is that the preconditioner is only a linear operator
if the coarse grained inverse is solved accurately.
In the numerical experiments described below we have found that the
coarse grained inverse only need be solved to an accuracy comparable to
that
desired for the overall solution on the fine grid.  This may be explained
by following heuristic argument.  Suppose the error from the inaccuracy in the
coarse grained inverse is random.  Then for a large problem,
dimension $O(10^6)$, the component of the error along a particular vector
is likely to be small.  Thus the error in orthogonality for the first few
iterations, $O(10)$, is small.  If only a small number of iterations are
required because of the quality of the approximate inverse then the
effect of small non-linearities of the preconditioner, introduced by
the recursive solvers, is negligible.

\section{Implementation}
To validate the algorithm we have implemented the solver in an object oriented
C++ based code.  The implementation is general enough to allow for an arbitrary
choice of solver and an arbitrary preconditioner at every level. We have tested
the code with the conjugate gradient solver using the three preconditioners
described in the previous section: polynomial preconditioner,
multi-grid preconditioner and the recursive multi-scale preconditioner.
In addition, the implementation includes the standard multi-grid algorithm.
This enabled us to compare the algorithms while using the same extension,
reduction, splitting and coarsening operators.

As a test case we used the porous flow equation
\begin{equation}
	\nabla \cdot (\K \cdot \nabla P) = S \;.
\end{equation}
Here, $\K$ is a permeability tensor, $P$  is the pressure and $S$ is a source.
We considered only a diagonal tensor and used the 
linear system derived from the standard
5-point stencil for the finite difference operator on a
two-dimensional regular grid.  In the discretization $P$ and $S$ are cell
centered fields whereas the components of $K$ are face centered;
from the cell centers $K_{xx}$ is offset by a $1/2$ cell in the x-direction
and $K_{yy}$ is offset by a $1/2$ cell in the y-direction.
The face centered components of the permeability field are obtained as
the harmonic mean of the adjacent cell centered values.
This discretization is a special case of support operator differencing
\cite{ShashkovSteinberg}
and results in a matrix that preserves
the positivity and symmetry of the differential operator.

The discrete operator has a unique decomposition, $\A = \D + \U + \L$
where $\D$ is diagonal, $\U$ is strictly upper triangular and
$\L$ is strictly lower triangular.  In terms of these matrices, we used
two operator splitting $\A = \P - \Q$:  the symmetric Gauss-Seidel splitting
for which $\P = (\D+ \L) \D^{-1} (\D+\U)$ and $\Q = \L \D^{-1} \U$,
and a modified Jacobi splitting with $\P =  2\D$ and $\Q = \D - \U - \L$.
In contrast to the standard Jacobi splitting, the modified version
is effective in damping short wavelength errors associated with
the checker board modes.  Though the Jacobi splitting does work with
our algorithm, it is not as effective and the results presented below
are for the Gauss-Seidel splitting.

We allow for both Dirichlet and Neumann boundary conditions.
Each point on the boundary can be independently chosen to satisfy one
or the other condition.  The type of boundary condition affects the
discretized $\K$ field along the boundary. 
In addition, the value of the boundary condition enters as
a source term in the boundary cells of the discretized equations.
The boundary source terms are large and under coarsening scale differently
then the source terms in the interior. To avoid complications from
the boundary source terms, we transform to a problem with zero
boundary conditions.  This is accomplished by generating a smooth
field $P_{\rm bf}$ that matches the boundary conditions.  We then
solve the problem for $\delta P = P - P_{\rm bf}$ with zero boundary
conditions but an additional source term
$-\nabla \cdot (\K \cdot \nabla P_{\rm bf})$.

The coarse grained permeability tensor is based on the transmissivities,
$T_{xx} = {\Delta y\over \Delta x} K_{xx}$ and
$T_{yy} = {\Delta x\over \Delta y} K_{yy}$ where $\Delta x$ and $\Delta y$
are the width and height of a grid cell.
Physically these are extensive rather than intensive quantities.
The transmissivity behaves like a conductance rather than a conductivity.
The $T_{xx}$ component is taken as
\begin{equation}
	{1 \over T_{xx}^c} = \sum_i {1\over \sum_j T_{xx}(i,j)}
\end{equation}
where the sum is over fine cells contained within a coarse cell.
When a fine cell only partially overlaps a coarse cell the value of $T_{xx}$
associated with the fine cell is scaled  by the fraction of the height
in the overlap and inversely as the fraction of the width in the overlap.
A similar construction is used for $T_{yy}$ (with the sums over i and j
interchanged).  The coarse grained operator remains diagonal.
It is typically not a scalar even if $\K$ on the finest grid is
chosen to be scalar.
This simple coarsening properly accounts for the effect of the boundary
conditions on the boundary $K$'s.  Furthermore, the form of the
Laplace operator is preserved under coarsening.

It is noteworthy that the scale factor between grids is not limited to
a factor of two, nor for that matter to an integer.  Consequently, the
number of grid points in a linear dimension does not have to be a power
of the scale factor.
Typically, we scaled between successive grid levels by a linear factor of 4
which reduces the dimension of the problem at each level by a factor of 16.
The effective scale factor between adjacent levels may
vary slightly in order to obtain an integer grid dimension.

The reduction operator is taken as the adjoint of the extension operator.
For the extension operator we use the tensor product of 1-D linear
interpolations.  A piecewise constant interpolation works nearly as well.
In contrast to the piecewise constant interpolation,
a linear interpolation requires a boundary condition.
For both Dirichlet and Neumann boundaries, we chose a zero slope for the
boundary interpolation.  While this correctly captures the Neumann
condition, it introduces a small error for the Dirichlet case.
However, the viability of the piecewise constant interpolation suggests
that the error is small and can be neglected.
Numerical experiments confirm this suggestion.
For a Neumann boundary, since $\delta P$ near the boundary can be large,
a zero boundary condition is highly detrimental to the convergence of the
algorithm.  The fact that the zero slope interpolation condition is acceptable
for both the Dirichlet and Neumann case greatly simplifies the
implementation of boundary conditions that vary between Dirichlet and Neumann
from cell to cell.

The degree $m$ of smoothing is related to the scale factor.
Since the smoothing operator $\Pinv \Q$ typically connects only neighboring
grid points we have chosen $m$ to be the same as the scale factor.
As a result the preconditioning matrix fully couples every fine grid cell.
Smaller values of $m$ would require a larger number of conjugate gradient
iterations for problems in which the short wavelengths dominate the solution.

Our implementation of the conjugate gradient algorithm is conventional,
as outlined in \cite[\S 2.3.1, Figure~2.5]{Templates}.  
We base the convergence criterion on the norm of the residual
rather than the norm with respect to the preconditioner.
This is because applying the preconditioner is the most expensive
operation in the conjugate gradient step.  It is performed at the
beginning of the cycle and hence is out of date when the check for convergence
is made at the end of the cycle.  A convergence criterion of the same type
is applied on every level.   The criterion for convergence on
the $k^{\hbox{th}}$ level is
\begin{equation}
	\frac{|| \vec r ||^2}{N_k} < f^k \times \epsilon^2
\end{equation}
where $\vec r$ is the residual, $N_k$ is the number of grid points,
$\epsilon$ is the desired root mean squared error of the residual on
the finest level, and $f$ is an adjustable parameter
that allows us to tighten the error criterion on the coarser levels.
The algorithm appears not to be very sensitive to the choice of $f$.
In practice we found $f = 0.1$ works well.  As $f$ increases the number
of iteration on the fine level gradually increases, while too small a
value of $f$ results in unnecessary iterations on the coarse levels.

\section{Numerical Examples}
We have tested our solver algorithm on several examples of the
porous flow equations.  The examples below
use both Dirichlet and Neumann boundary conditions, typically,
constant pressure $P=1$ on the left and $P=0$ on the right, and no flow
on the top and bottom.  We generated random log-normal permeability field with
either a Gaussian auto-correlation function or a power law auto-correlation
function
\[
	\left({1\over 1+ (\vec r_1-\vec r_2) \mathbf{\Lambda}
				(\vec r_1-\vec r_2)}\right)^{1\over 4}
\]
where $\mathbf{\Lambda}$ is a positive definite matrix which defines
cutoff lengths for the power law behavior.
The variance of the permeability field is adjusted by scaling the log
of the field.  Similar problems have previously been used to test the conjugate
gradient algorithm with a multi-grid preconditioner \cite{AshbyFalgout}.

Our test examples include meshes varying in size by a factor of $1000$;
from $2\times10^3$ to $2\times 10^6$ grid points.
We also have varied the difficulty of the problem by
increasing the variance of the permeability field to obtain a maximum
to minimum permeability ratio up to $8\times10^9$.
In addition, we have adjusted the error tolerance to vary the
accuracy of the solution up to machine accuracy.

\subsection{Base case}
\begin{figure}[t] 
\centerline{\psfig{file=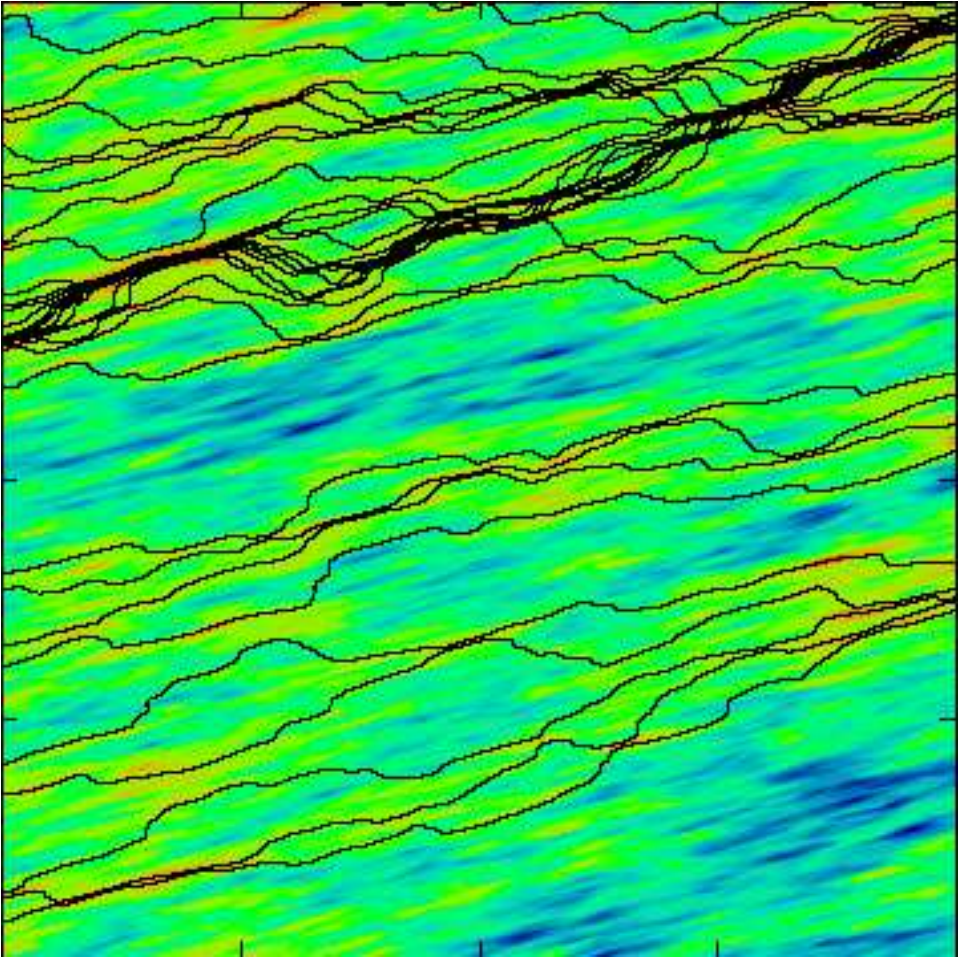,height=5truein}}
\caption{
	Flow lines (in black) superimposed on the log of permeability field
	for the base case.
	The permeability field is random log-normal with zero mean, a
	variance of 2 and a power law auto-correlation.  The field is
	discretized on a $1000\times 1000$ grid and its log ranges
	from a low of -7.3 (blue) to a high of 7.9 (red).
	We note for this ``fractal'' permeability field the flow lines
	tend to follow in narrow channels.  This is in contrast to our
	experience with Gaussian auto-correlations.
	(The raggedness of the flow lines is a consequence of the
	compression algorithm used to reduce the size of the plot file
	and is not an indication of either the resolution or
	the accuracy of the solution.)
}							\label{fig:FlowLines} 
\end{figure}

For a basic test case we used a power law field on a quarter of the unit square.
We choose a variance of 2 with a zero mean and minimum correlation lengths
of 0.016 and 0.002 oriented at 15 degrees with respect to the x-axis.
The resulting permeability field contains a wide range of wavelengths.
On a $1000\times1000$ grid the distribution of the log of the
permeability field extends over 3.5 standard deviations.  Consequently,
the permeability field varies from~$10^{-3}$ to~$10^3$ and the ratio of
its maximum to minimum value is $10^6$.
The flow lines for the solution superimposed on the log of the
permeability field are shown in
figure~\ref{fig:FlowLines}.
The conjugate gradient solver on the finest level required 5 iterations
to reduce the mean squared residual by a factor of $10^{10}$.
Later we show that this corresponds to 5 digits accuracy.
Replacing the Neumann boundary conditions on the top and bottom
with Dirichlet boundary conditions does not change the performance
of the solver.

\begin{figure}[t] 
\centerline{\psfig{file=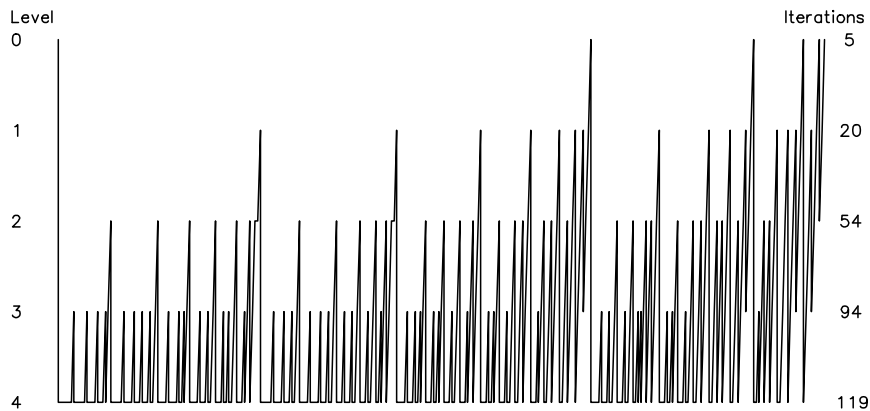,width=5truein}}
\caption{
	Level tree showing the dynamic transition between levels
	for the base case.  Level 0 is the finest grid and
	level 4 is the coarsest grid.
}							\label{fig:LevelTree} 
\end{figure}

The algorithm dynamically determines its effort on every level.
From the level tree shown in
figure~\ref{fig:LevelTree}
it is seen that the effort on the coarsest levels diminishes as the calculation
proceeds.  The total number of iterations is summarized in table~\ref{table1}.
The effort on each level is proportional to the total number of iterations
on that level times the dimension of the level. It is seen from the
table that the total effort is dominated by the computation on the finest level.

To set the scale, the total computational effort is equivalent to about 1600
scalar products of vectors on the finest grid.
On a workstation (SUN Ultra~I, 170~Mhz) or a PC (PentiumPro, 200~Mhz)
the time per point is $240\,\mu$s.
This time is for the solver only and does not include initialization of
the permeability fields on the coarse grids.
Our implementation is memory efficient.
The solver requires storage for a total of 10 vectors:
3 for the matrix, one each for the pressure field, the source field,
the residual and the conjugate direction, and three temporaries to 
evaluate the preconditioner (\ref{eq:2.10}) on the residual.

\begin{table}
\caption{Iteration count on each level for the base case.
	Computational effort on a level is proportional
	to the dimension of the level times the total number of
	iterations on that level.
	}							\label{table1} 

	\begin{center} \footnotesize
	\begin{tabular}{r| r @{ $\times$ } r @{ = }r|r|r|r} \hline 
	Level & \multicolumn{2}{|c}{Grid} &dimension &iterations
			& iterations$\times$ dimension & per cent total\\ \hline
	   4  &   4 &   4 &     16	&119	&1904		&0.03  \\
	   3  &  16 &  16 &    256	&94	&24064		&0.37  \\
	   2  &  63 &  63 &   3869 	&63	&214326		&3.29  \\
	   1  & 252 & 252 &  63504	&20	&1270080	&19.48 \\
	   0  &1001 &1001 &1002001	&5	&5010005	&76.84 \\ \hline
	    \multicolumn{6}{r}{total 6520379}	\\
	\end{tabular}
	\end{center} 
\end{table}

Specializing the algorithm to the
Laplace equation reduces the computational effort to $50\,\mu$s per point
and decreases the required memory by 4 vectors.
The difference in speed, as well as memory, can be attributed largely
to the much more efficient coding for the simpler linear operator.  It
is a remarkable fact that the same preconditioned conjugate gradient
algorithm works equally well for an operator with a rapidly varying
permeability field as for the simple Laplace operator.

We choose this problem, which is fairly large for a workstation,
as a base case because we felt it is necessary to have several levels
of coarsening to assess the quality of the algorithm.
With a scale factor of roughly 4,
there are only 4 coarsenings for a $1000\times 1000$ grid.
The problem we have choosen is not trivial.  The polynomial preconditioner
fails to converge, and so does the standard multigrid algorithm.  Our
implementation of Tatebe's algorithm succeeds and will be discussed in
more detail later.

\begin{figure}[t] 
\centerline{\psfig{file=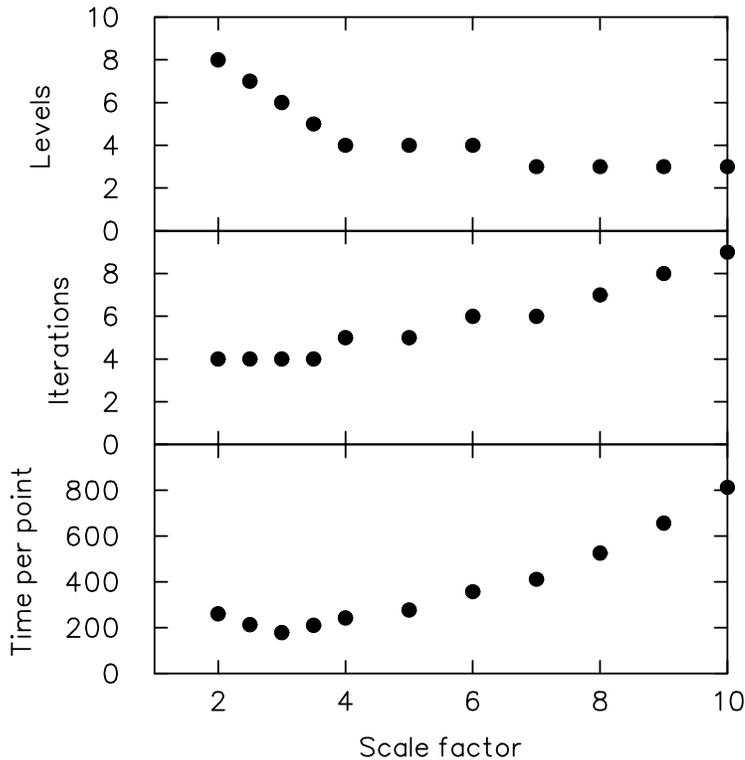,width=4.truein}}
\caption{
	The effect of the scale factor between grids on the base case.
	Shown are the number of levels, iterations on the finest grid,
	and the time per point.
}							\label{fig:ScaleFactor} 
\end{figure}

Though the algorithm has not been optimized, it is relatively insensitive
to the parameters characterizing the algorithm.  Figure~\ref{fig:ScaleFactor}
shows the effect of the scale factor between grids.  For scale factors
between 2 and 5, the time per point varies by only 30\%.
A minimum time of $180\,\mu$s per point occurs with a scale factor of 3.

\subsection{Scaling behavior} 
For high resolution, the scaling properties of the solver algorithm
are critical.  To test the scaling behavior we solved a series of problems
on sucessively larger grids varying in dimension by 3 orders of magnitude.
We found the effort per point to be independent of the size of the grid.

We started by generating a random permeability field on a unit square
with a very fine grid, $2048 \times 2048$.  Then the grid is truncated
to a square subgrid by using only the subregion defined by the lower left
corner and a specified upper right corner.
The $1000\times1000$ field corresponds to our base case
and is choosen to have a variance of 2 with a mean of zero.
The cutoff lengths of the correlation function on this sequence of
grids remains constant
(32 by 4 cells) but the variance of the permeability field increases with
the size of the grid.  Consequently, the problem gets harder as the
size increases.

\begin{figure}[t] 
\centerline{\psfig{file=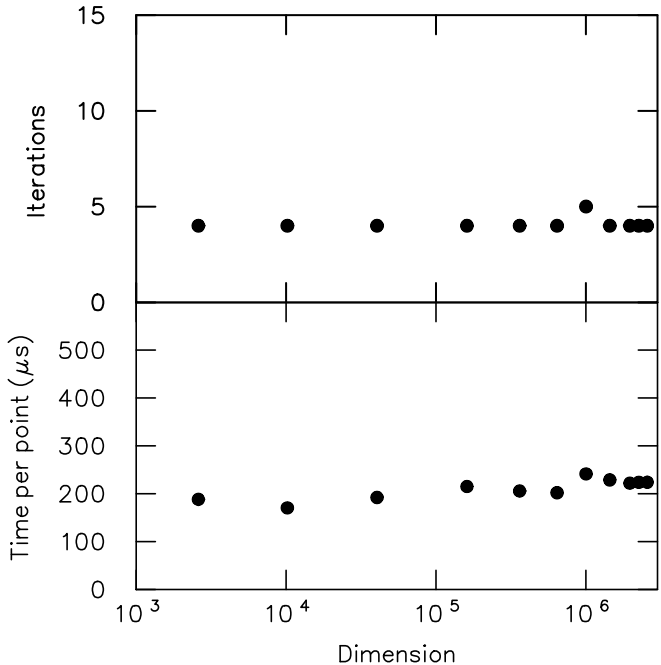,height=3.5truein}}
\caption{
	Number of iterations on the finest grid and time per point
	versus number of grid points for the recursive multi-scale
	conjugate gradient algorithm.
}							\label{fig:scaling} 
\end{figure}

Figure~\ref{fig:scaling} shows the scaling behavior on grids ranging
from $50\times 50$ to $1600\times 1600$.  It can be seen that
the number of iterations and the time per point are almost constant.
This indicates that the computational effort scales linearly with the
grid dimension.  Despite the variable coefficients, this scaling is better
than the $N \log N$ scaling for solving the Laplace equation with fast
Fourier transforms.

\begin{figure}[t] 
\centerline{\psfig{file=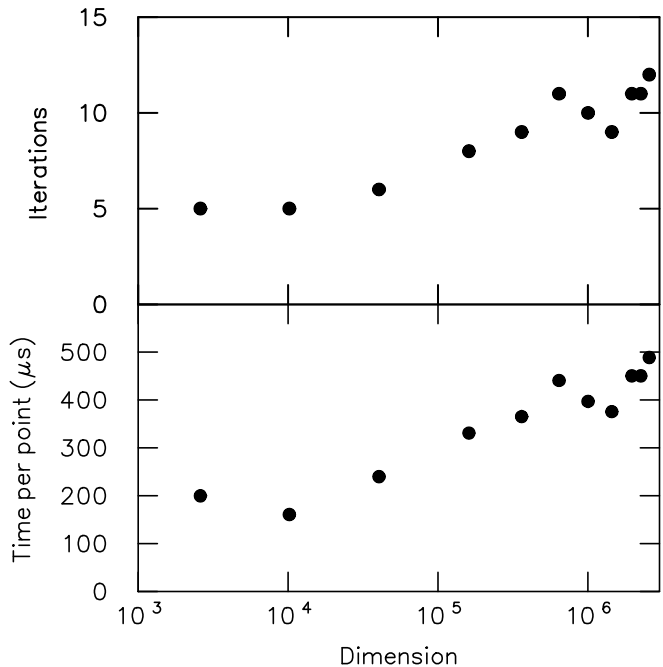,height=3.5truein}}
\caption{
	Number of iterations on the finest grid and time per point
	versus number of grid points for the multi-grid preconditioned
	conjugate gradient (Tatebe's) algorithm.
}						\label{fig:scaling_Tatebe} 
\end{figure}

We have run the same problems with Tatebe's multi-grid preconditioning
algorithm.  As seen in figure~\ref{fig:scaling_Tatebe} this algorithm
requires more iterations. The overall trend is linear which indicates
that the computational effort scales as $N\log N$. 
This is in line with the results of Ashby \& Falgout \cite{AshbyFalgout}.
A comparison with their results is necessarily imprecise because
they concentrated on 3-D problems and a parallel implementation.
The larger number of iteration in Tatebe's algorithm is partially offset
by a lower cost per iteration.  For the largest problems in this series,
the time per point is a factor of two greater than for our algorithm.
For the small problems the two algorithms are comparable in time.

In addition, we have tested the scaling behavior on a family of problems
with the permeability field generated by interpolating
from the very fine grid to coarser grids on the same
physical domain and then rescaling to obtain the same variance.
For this family of problems, the cutoff length to the auto correlation function
in units of cells
is proportional to the grid size.  The decreasing smoothness of the
discretized field, in terms of cell to cell variation, increases the difficulty of
the problem for the smaller grids.  Due to the changing correlation
length, with this set of problems the time per point for our algorithm
decreases slightly as the grid dimension increases while Tatebe's algorithm
continues to scale as $N\log N$.  For the largest grid our algorithm
is again 50\% faster than Tatebe's algorithm.

The numerical evidence suggest that our algorithm scales as ${\cal O}(N)$.
This linear scaling behavior is one of the strong points of the
recursive multi-scale conjugate gradient algorithm.
In fact the only algorithms robust enough to solve this class of
problems are the hybrids that combine multi-scaling with conjugate gradient.
The simplest representative of this hybrid class is Tatebe's algorithm.
Our algorithm is an example of the trade-off between a higher cost per
iteration and a more accurate preconditioner.
These examples show that this tradeoff reduces the computational time
and also improve the scaling behavior compared to Tatebe's algorithm.
Because of the favorable scaling, we expect that the advantage
of the recursive multi-scale conjugate gradient algorithm to increase
with the problem size.

\subsection{Accuracy} 
In order to obtain an estimate of the accuracy of our method,
we first generated a test problem for which the exact solution
is known.  To this end we solved our base case approximately
using another solver.  We intentionally did not strive for high accuracy.
By adding the small residual of the approximate solution to the source term
of the base problem, we created a new test problem which by construction
is solved exactly by the approximate pressure field of the base problem.

\begin{figure}[t] 
\centerline{\psfig{file=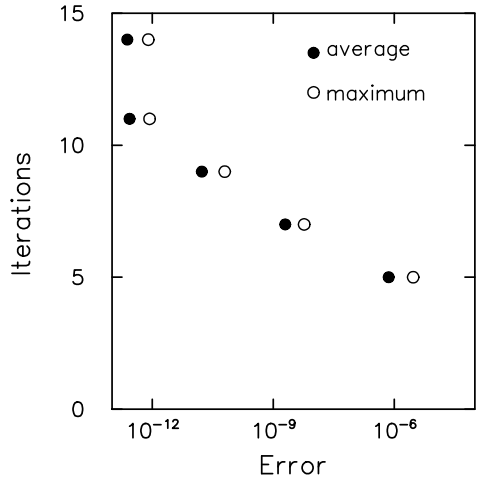,height=2.5truein}}
\caption{Iterations vs accuracy for the base case.  The solution field is
of order 1.
}						\label{fig:accuracy} 
\end{figure}

By varying the error tolerance of our solver we generated a sequence of
solutions.  For these solutions, the iteration count as a function of
the accuracy is shown in figure~\ref{fig:accuracy}.
Both the root mean squared error and the maximum error
are used as measures of the accuracy.
We observe that the iteration
count increases linearly with the log of the accuracy until the improvement
of the solution is limited by machine accuracy.
The linear trend indicates that
each iteration reduces the error by a factor of about 11.
Consequently, only a small number of iterations are needed
to achieve machine accuracy.
The fact that the maximum error is within a factor of 3 of the root mean squared
error indicates that the solution is uniformly accurate.  This is another
strong point of our algorithm and is a consequence
of the preconditioner acting on all length scales.

\begin{figure}[t] 
\centerline{\psfig{file=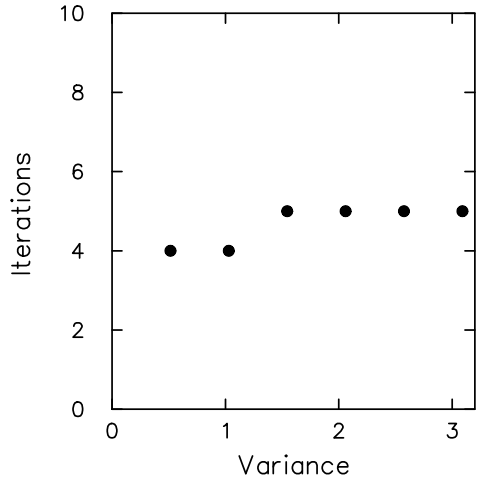,height=2.5truein}}
\caption{
Iterations vs variance of the log of the permeability field.
The variance is adjusted by scaling the field of the base case.
}						\label{fig:variance} 
\end{figure}

To further test the robustness of our algorithm we increased the difficulty
of the base problem by increasing the variance of the permeability field.
The iteration count as a function of the variance is shown
in figure~\ref{fig:variance}.  We observe that the iteration count is
almost constant.  A variance of 3 results in the value of permeability
field varying from minimum to maximum by a factor of $8 \times 10^9$.
For the larger variances, round-off errors in the finite difference
approximation to $\nabla \cdot (\K \cdot \nabla P)$ limits the accuracy
to which the solution can be computed.

These tests shows that our algoithm is very robust.  The preconditioner
is effective on difficult problems.
With a small increase in the number of iterations the solution can be
driven to near machine accuracy.
It is remarkable, that despite rather large spatial variation in the value of
the permeamility field, a solution can still be obtained in only 5 iterations.

\subsection{Non-isotropic permeability field} 
We performed preliminary tests to determine whether the algorithm
can be applied to problems with large aspect ratios in either the domain
or the individual grid cells.  
For the permeability field we used a strip of the fine field generated for
the scaling study (full width and 25\% of the height).
The left half of the new field
corresponds to the top half of the permeability field of the base case.
The discretized grid of $2000\times 500$ has the same overall dimension,
$10^6$, as the base case.
In addition, as with the base case, the variance is set to~2.

We then reinterpreted the discretized field by assuming a
cell aspect ratio of~10 to~1.  While the initial field ranged from 0 to 1
in both the $x$- and $y$-directions, the new field ranges from 0 to 10 in
the x-direction and from~0.25 to~0.5 in the y-direction.
The rescaling results in a new
permeability field in a channel with a length 40 times its width and
in which the typical features (ratio of the correlation lengths)
have an aspect ratio of 80 to 1.

\begin{figure}[t] 
\centerline{\psfig{file=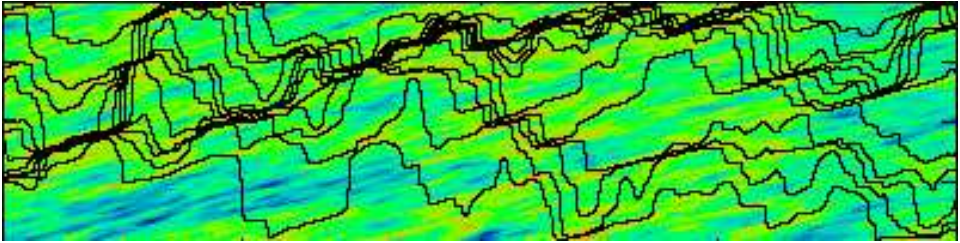,width=5truein}}
\caption{
	Flow lines (in black) superimposed on the log of permeability field
	for the flow in a narrow channel.
	The permeability field is random log-normal with zero mean, a
	variance of 2 and a power law auto-correlation.  The field is
	discretized on a $2000\times 500$ grid and its log ranges
	from a low of -7.2 (blue) to a high of 7.9 (red).
	The cells have a 10 to 1 aspect ratio resulting in a domain
	with a 40 to 1 aspect ratio.
	To fit on the page, the x direction has been
	compressed by a factor of 10.  This is equivalent to a non-isotropic
	permeability field with $K_{yy}$ a factor of 100 times $K_{xx}$.
	(The raggedness of the flow lines is a consequence of the
	compression algorithm used to reduce the size of the plot file
	and is not an indication of either the resolution or
	the accuracy of the solution.)
	}						\label{fig:Strip} 
\end{figure}

The matrix for the problem with a 10 to 1 cell aspect ratio is equivalent
to the matrix corresponding to square cells and a non-isotropic
permeability field with $K_{yy}$ a factor of 100 times $K_{xx}$.
The flow lines superimposed on the log of the permeabilty field
are shown in figure~\ref{fig:Strip}.  The effect of the anisotropy
is seen in the abrupt changes in direction of the flow lines in order
to follow paths of high permeability.

We found that the solver performed better when the coarsening strategy
aimed for cells with an aspect ratio of 1.  Thus, despite the smaller
dimension, the grid is first coarsened in the y-direction, and then
uniformly in both directions.  This is similar to the
``semi-coarsening'' strategy used for multi-grid algorithms.
The iteration count on each level is shown in table~\ref{table2}.
The finest level still required only 4 iterations but the semi-coarsening
does not reduce the grid dimension by as large a factor for the first
two levels and
the time per point of 768$\,\mu$s is 2.4 times as large as for the base
case.  We conjecture that a better coarsening algorithm which
accounts for off-diagonal components of the permeability tensor
would not require the semi-coarsening and that with a uniform reduction
in grid size the algorithm would be as efficient as for the base case.

\begin{table}
\caption{Iteration count on each level for the flow in a narrow channel.
	Computational effort on a level is proportional
	to the dimension of the level times the total number of
	iterations on that level.
	}							\label{table2} 

	\begin{center} \footnotesize
	\begin{tabular}{r| r @{ $\times$ } r @{ = }r|r|r|r} \hline 
	Level & \multicolumn{2}{|c}{Grid} &dimension &iterations
			& iterations$\times$ dimension & per cent total\\ \hline
	   4  & 160 &  4 &    640	&1548	&990720	 & 6.3	\\
	   3  & 520 & 13 &   6760	& 191	&1291160 & 8.2	\\
	   2  &1760 & 44 &  77440 	&  66	&5111040 &32.3	\\
	   1  &2000 &147 & 294000	&  15	&4410000 &27.9	\\
	   0  &2000 &500 &1000000	&   4	&4000000 &25.3	\\ \hline 
	    \multicolumn{6}{r}{total 15802920}	\\
	\end{tabular}
	\end{center} 
\end{table}

\section{Conclusion}
For a large system, any iterative algorithm requires an excellent
preconditioner.  For PDEs that can be reasonably described by discretization,
a coarser discretization forms the basis for a good preconditioner.
If our algorithm is well suited for the finest level of discretization,
it stands to reason
that it is equally well suited to the next level.  Thus we are naturally
led to a class of preconditioned algorithms which are recursive in nature.
The preconditioner applies the very same algorithm on a coarser scale 
until the problem is either so coarsely resolved that further coarsening
is detrimental or the problem is sufficiently small for direct solvers
to be more efficient.

We have applied this philosophy to the conjugate gradient algorithm.
However, it is clear that the same ideas apply to other solvers.
Indeed, the high quality of the preconditioner is even more valuable
for other Krylov space methods which need to store a set of conjugate
directions.  In our example the  memory required is equivalent to 10
solution vectors on the finest grid.
Since the algorithm generally converged in fewer than 10 iterations, the
overhead of storing intermediate vectors would increase the memory
requirement by less than a factor of~2.
This suggests that similar algorithms would be effective when
the operator is not symmetric or not positive definite,
significantly broadening the
number of problems that can be addressed.  Even for porous flow
problems of the type addressed here, optimum coarsening techniques
would naturally lead to nonsymmetric  $\K$ tensors.

Our solver is yet another example that demonstrates the advantages of
hybrid schemes that combine multi-grid and conjugate gradient algorithms.
Because of their scaling
behavior and robustness, hybrids algorithms are the only ones effective
on truly large problems.  Our implementation allows us to define a
different preconditioner on every level.  The algorithm can be
specialized to a particular problem by tailoring the solver on every
level to the frequency distribution of scales.  From this point of view
Tatebe's algorithm and our algorithm are the end points of a large
class of multi-scale algorithms.
Our algorithm has the advantage of a favorable scaling,
but for some problems Tatebe's algorithm may be more suitable.

The ideas on which our approach is based are quite general and
transcend the specific implementation.   In many regards our implementation
is a particularly simple example of the general approach. We expect in
the future to see similar recursive algorithms that have equally
favorable scaling behavior but can be applied to a large class of
problems with more complex  gridding and coarse graining schemes.
In this paper we were also not concerned with  parallelizing.
However, it is clear that similar approaches can be parallelized
and that a significant increase in speed could thus be obtained.

\bibliography{Refs}

\end{document}